\documentclass[preprint,12pt]{elsarticle}
\usepackage{amsmath}

\usepackage{amssymb}
\begin{document}

\begin{frontmatter}


\ead{aaalikhanov@gmail.com}

\title{Numerical methods of solutions of
boundary value problems for the multi-term variable-distributed
order diffusion equation}


\author{Anatoly A. Alikhanov}

\address{Kabardino-Balkarian State University, ul. Chernyshevskogo 173,
 Nalchik,  360004,   Russia}

\begin{abstract}
Solutions of the Dirichlet and Robin boundary value problems for the
multi-term variable-distributed order diffusion equation are
studied. A priori estimates for the corresponding differential and
difference problems are obtained by using the method of the energy
inequalities. The stability and convergence of the difference
schemes follow from these a priory estimates. The credibility of the
obtained results is verified by performing numerical calculations
for test problems.

\end{abstract}


\end{frontmatter}

% \linenumbers


%%%%%%%%%%%%%%%%%%%%%%%
\section{Introduction}
%%%%%%%%%%%%%%%%%%%%%%

Differential equations with fractional order derivatives provide a
powerful mathematical tool for accurate and realistic description of
physical and chemical processes proceeding in media with fractal
geometry \cite{Nakh:03,Podlub:99,Hilfer:00,Kilbas:06,Uchaikin:08}.
It is known that the order of a fractional derivative depends on the
fractal dimension of medium \cite{Kobelev:98,Kobelev:99}. It is
therefore reasonable to construct mathematical models based on
partial differential equations with the variable and distributed
order derivatives
\cite{Nakh:03,Atanack:11,Coimb:03,Lorenzo:02,Podlubny:12,Luchko:09,Pskh:05,Pskh:04}.
Analytical methods for solving  such equations are scarcely
effective, so that the development of the corresponding numerical
methods is very important.

The initial-boundary-value problems for the generalized multi-term
time fractional diffusion equation over an open bounded domain
$G\times(0,T), G\in \Bbb{R}^n$ were considered \cite{Luchko_mt}.
Multi-term linear and non-linear diffusion-wave equations of
fractional order were solved in \cite{pdu_1} using the Adomian
decomposition method.  Applications of  the homotopy analysis and
new modified homotopy perturbation methods to solutions of
multi-term linear and nonlinear diffusion-wave equations of
fractional order are discussed in \cite{Jafari_1,Jafari_2}. In the
papers \cite{liu_1,liu_2} analytical solutions for the multi-term
time-fractional diffusion-wave and the multi-term time-space
Caputo-Riesz fractional advection-diffusion equations on a finite
domain are studied. The fundamental solution of the multi-term
diffusion equation with the Dzharbashyan-Nersesyan fractional
differentiation operator with respect to the time variables is
constructed in \cite{Pskhu_pdu}.

Several methods for solving variable and distributed order
fractional differential equations with various kinds of the variable
and distributed fractional derivative have been proposed
\cite{Podlubny:12,Liu:12_1,Liu:12_2,Liu:10,Liu:09,Liu:09_2,Diethelm:09,Stojanovic:11}.
 A priory estimates for the difference problems obtained in
 \cite{ShkhTau:06,ShkhLaf:09,ShkhLaf:10} by using the maximum
 principle imply the stability and convergence of the considered
 difference schemes.
  Using the energy inequality method, a priori estimates for the
solution of the Dirichlet and Robin boundary value problems for the
fractional and variable order diffusion equation with Caputo
fractional derivative have been obtained \cite{Alikh:10},
\cite{Alikh:12}.

%%%%%%%%%%%%%%%%%%%%%%%%%%%%%%%%%%%%%%%%%%%%%%%%%%%%%%%%
\section{Boundary value problems in differential setting}
%%%%%%%%%%%%%%%%%%%%%%%%%%%%%%%%%%%%%%%%%%%%%%%%%%%%%%%%%

%%%%%%%%%%%%%%%%%%%%%%%%%%%%%%%%%%%%%%%%%%%%%%%%%%%%%
\subsection{The Dirichlet boundary value problem}
%%%%%%%%%%%%%%%%%%%%%%%%%%%%%%%%%%%%%%%%%%%%%%%%%%%%

In rectangle $\bar Q_T=\{(x,t): 0\leq x\leq l, 0\leq t\leq T\}$ let
us study the boundary value problem
\begin{equation}\label{ur1}
\Bbb{P}_{
(\omega)}^{(\theta)}\left(\partial_{0t}\right)u(x,t)=\frac{\partial
}{\partial x}\left(k(x,t)\frac{\partial u}{\partial
x}\right)-q(x,t)u+f(x,t) ,\,\, 0<x<l,\,\, 0<t\leq T,
\end{equation}
\begin{equation}
u(0,t)=0,\quad
 u(l,t)=0,\quad 0\leq t\leq T, \label{ur2}
\end{equation}
\begin{equation}
u(x,0)=u_0(x),\quad 0\leq x\leq l, \label{ur3}
\end{equation}
where
$$
\Bbb{P}_{
(\omega)}^{(\theta)}\left(\partial_{0t}\right)u(x,t)=\int\limits_{\alpha}^{\beta}d\gamma\sum\limits_{r=1}^{m}\omega_r(x,\gamma)\partial_{0t}^{\theta_r(x,\gamma)}u(x,t),
$$
$$ \alpha<\beta, \quad 0<\theta_r(x,\gamma)<1, \quad
\omega_r(x,\gamma)\geq 0, \quad
 r=1,2,...,m, \quad \text{for all}
$$
$$  (x,\gamma)\in[0,l]\times[\alpha,\beta],
\quad
\int\limits_{\alpha}^{\beta}d\gamma\sum\limits_{r=1}^{m}\omega_r(x,\gamma)>0,\quad
\theta_r(x,\gamma)\in C[0,l]\times[\alpha,\beta],
$$
$$
0<c_1\leq k(x,t)\leq c_2,\quad q(x,t)\geq 0,
$$
$\partial_{0t}^{\theta_r(x,\gamma)}u(x,\eta)=\int_{0}^{t}u_{\eta}(x,\eta)(t-\eta)^{-\theta_r(x,\gamma)}d\eta/\Gamma(1-\theta_r(x,\gamma))$
 is a  Caputo fractional derivative of order $\theta_r(x,\gamma)$ \cite{Cap:69,MainGor:07}.

 The existence of the solution for the initial boundary value
problem of fractional, multi-term and distributed order diffusion
equation has been proven in
\cite{Luchko:09,liu_1,Luchko:10,Luchko:11,Meersch:09,liu_111}.

Let us assume further the existence of a solution  $u(x,t)\in
C^{2,1}(\bar Q_T)$ for the problems  (\ref{ur1})--(\ref{ur3}), where
$C^{m,n}$  is the class of functions, continuous together with their
partial derivatives of the order  $m$ with respect to  $x$ and order
$n$ with respect to  $t$ on $\bar Q_T$.

{\bf Lemma 1.} For any functions $v(t)$ and $w(t)$ absolutely
continuous on $[0,T]$, one has the equality:

$$
v(t)\Bbb{P}_{
(\bar\omega)}^{(\bar\theta)}\left(\partial_{0t}\right)w(t)+
w(t)\Bbb{P}_{
(\bar\omega)}^{(\bar\theta)}\left(\partial_{0t}\right)v(t)=\Bbb{P}_{
(\bar\omega)}^{(\bar\theta)}\left(\partial_{0t}\right)(v(t)w(t))+
$$

\begin{equation}
+\int\limits_{\alpha}^{\beta}d\gamma\sum\limits_{r=1}^{m}\frac{\bar\omega_r(\gamma)\bar\theta_r(\gamma)}{\Gamma(1-\bar\theta_r(\gamma))}\int\limits_{0}^{t}\frac{d\xi}{(t-\xi)^{1-\bar\theta_r(\gamma)}}
\int\limits_{0}^{\xi}\frac{v'(\eta)d\eta}{(t-\eta)^{\bar\theta_r(\gamma)}}\int\limits_{0}^{\xi}\frac{w'(s)ds}{(t-s)^{\bar\theta_r(\gamma)}},
 \label{ur4}
\end{equation}
where $\bar\omega_r(\gamma)\geq0$, $0<\bar\theta_r(\gamma)<1$,  for
all $\gamma\in [\alpha,\beta]$,
$\int_{\alpha}^{\beta}d\gamma\sum\limits_{r=1}^{m}\bar\omega_r(\gamma)>0$.

{\bf Proof.}  For any fixed $\gamma\in[\alpha,\beta]$ and
$r\in\{1,2,\ldots,m\}$, relying on lemma 1~\cite{Alikh:12} one finds
the following equality
$$
v(t)\partial_{0t}^{\bar\theta_r(\gamma)}w(t)+
w(t)\partial_{0t}^{\bar\theta_r(\gamma)}v(t)=\partial_{0t}^{\bar\theta_r(\gamma)}(v(t)w(t))+
$$

\begin{equation}
+\frac{\bar\theta_r(\gamma)}{\Gamma(1-\bar\theta_r(\gamma))}\int\limits_{0}^{t}\frac{d\xi}{(t-\xi)^{1-\bar\theta_r(\gamma)}}
\int\limits_{0}^{\xi}\frac{v'(\eta)d\eta}{(t-\eta)^{\bar\theta_r(\gamma)}}\int\limits_{0}^{\xi}\frac{w'(s)ds}{(t-s)^{\bar\theta_r(\gamma)}}.
 \label{ur4.11}
\end{equation}

 Multiplying (\ref{ur4.11}) by $\bar\omega_r(\gamma)$ and summing
 the result over $r$ from $1$ to $m$, then integrating over $\gamma$ from $\alpha$ to $\beta$ one
obtains (\ref{ur4}). The proof of the lemma 1 is complete.

{\bf Corollary.} For any function $v(t)$ absolutely continuous on
$[0,T]$, the following equality takes place:
$$
v(t)\Bbb{P}_{
(\bar\omega)}^{(\bar\theta)}\left(\partial_{0t}\right)v(t)=
\frac{1}{2}\Bbb{P}_{
(\bar\omega)}^{(\bar\theta)}\left(\partial_{0t}\right)v^2(t)+
$$
\begin{equation}
+\int\limits_{\alpha}^{\beta}d\gamma\sum\limits_{r=1}^{m}\frac{\bar\omega_r(\gamma)\bar\theta_r(\gamma)}{2\Gamma(1-\bar\theta_r(\gamma))}\int\limits_{0}^{t}\frac{d\xi}{(t-\xi)^{1-\bar\theta_r(\gamma)}}
\left(\int\limits_{0}^{\xi}\frac{v'(\eta)d\eta}{(t-\eta)^{\bar\theta_r(\gamma)}}\right)^2,
 \label{ur4.200}
\end{equation}
where $\bar\omega_r(\gamma)\geq0$, $0<\bar\theta_r(\gamma)<1$,  for
all $\gamma\in [\alpha,\beta]$,
$\int_{\alpha}^{\beta}d\gamma\sum\limits_{r=1}^{m}\bar\omega_r(\gamma)>0$.

Let us use the following notation:
$\|u\|_0^2=\int\limits_{0}^{l}u^2(x,t)dx$,
$D_{0t}^{-\nu}u(x,t)=\int\limits_{0}^{t}(t-s)^{\nu-1}u(x,s)ds/\Gamma(\nu)$
is a fractional Riemann-Liouville integral of order  $\nu>0$.

 {\bf Theorem 1.} If $k(x,t)\in C^{1,0}(\bar{Q}_T)$, $q(x,t), \,
f(x,t)\in C(\bar Q_T)$, $k(x,t)\geq c_1>0$, $q(x,t)\geq 0$
everywhere on $\bar Q_T$, then the solution  $u(x,t)$ of the problem
(\ref{ur1})--(\ref{ur3}) satisfies the a priori estimate:
$$
\int\limits_{0}^{l}\Bbb{P}_{
(\omega)}^{(\theta-1)}\left(D_{0t}\right)u^2(x,t)dx+
c_1\int\limits_{0}^{t}\|u_x(x,s)\|_0^2ds\leq
$$
\begin{equation}\label{ur4.300}
\leq \frac{l^2}{2c_1}\int\limits_{0}^{t}\|f(x,s)\|_0^2ds+
\int\limits_{0}^{l}u_0^2(x)dx\int\limits_{\alpha}^{\beta}d\gamma\sum\limits_{r=1}^{m}\frac{\omega_r(x,\gamma)t^{1-\theta_r(x,\gamma)}}{\Gamma(2-\theta_r(x,\gamma))},
\end{equation}
where
$\Bbb{P}_{(\omega)}^{(\theta-1)}\left(D_{0t}\right)=\int_{\alpha}^{\beta}d\gamma\sum\limits_{r=1}^{m}\omega_r(x,\gamma)D_{0t}^{\theta_r(x,\gamma)-1}$.

{\bf  Proof.} Let us multiply equation (\ref{ur1}) by $u(x,t)$ and
integrate the resulting relation over $x$ from $0$ to $l$:
$$
\int\limits_{0}^{l}u(x,t)\Bbb{P}_{(\omega)}^{(\theta)}\left(\partial_{0t}\right)u(x,t)dx
-\int\limits_{0}^{l}u(x,t)(k(x,t)u_x(x,t))_xdx+
$$
\begin{equation}\label{ur5}
+\int\limits_{0}^{l}q(x,t)u^2(x,t)dx=\int\limits_{0}^{l}u(x,t)f(x,t)dx.
\end{equation}

Then transform the terms in identity (\ref{ur5}) as
\begin{equation}\label{ur6}
-\int\limits_{0}^{l}u(x,t)(k(x,t)u_x(x,t))_xdx=\int\limits_{0}^{l}k(x,t)u_x^2(x,t)dx\geq
c_1\|u_x(x,t)\|_0^2,
\end{equation}
\begin{equation}\label{ur7}
\left|\int\limits_{0}^{l}u(x,t)f(x,t)dx\right|\leq
\varepsilon\|u(x,t)\|_0^2+\frac{1}{4\varepsilon}\|f(x,t)\|_0^2,
\quad \varepsilon>0.
\end{equation}

Using the equality (\ref{ur4.200}) one obtains

\begin{equation}\label{ur9}
\int\limits_{0}^{l}u(x,t)\Bbb{P}_{(\omega)}^{(\theta)}\left(\partial_{0t}\right)u(x,t)dx\geq
\frac{1}{2}\int\limits_{0}^{l}\Bbb{P}_{(\omega)}^{(\theta)}\left(\partial_{0t}\right)u^2(x,t)dx.
\end{equation}

Taking into account the above performed transformations, from the
identity (\ref{ur5}) one arrives at the inequality
\begin{equation}\label{ur10.5}
\frac{1}{2}\int\limits_{0}^{l}\Bbb{P}_{(\omega)}^{(\theta)}\left(\partial_{0t}\right)u^2(x,t)dx+
c_1\|u_x(x,t)\|_0^2\leq\varepsilon\|u(x,t)\|_0^2+\frac{1}{4\varepsilon}\|f(x,t)\|_0^2.
\end{equation}

Using the inequality  $\|u(x,t)\|_0^2\leq (l^2/2)\|u_x(x,t)\|_0^2$,
from the inequality (\ref{ur10.5}) at  $\varepsilon=c_1/l^2$, one
obtains
\begin{equation}\label{ur11}
\int\limits_{0}^{l}\Bbb{P}_{(\omega)}^{(\theta)}\left(\partial_{0t}\right)u^2(x,t)dx+
c_1\|u_x(x,t)\|_0^2\leq\frac{l^2}{2c_1}\|f(x,t)\|_0^2.
\end{equation}

Changing the variable $t$ by $s$ in inequality (\ref{ur11}) and
integrating it over $s$ from $0$ to $t$, one obtains the a priori
estimate (\ref{ur4.300}).

The uniqueness and the continuous dependence of the solution of the
problem (\ref{ur1})--(\ref{ur3}) on the input data follow from the a
priori estimate (\ref{ur4.300}).

%%%%%%%%%%%%%%%%%%%%%%%%%%%%%%%%%%%%%%%%%%%%%%%%
\subsection{The Robin boundary value problem.}
%%%%%%%%%%%%%%%%%%%%%%%%%%%%%%%%%%%%%%%%%%%%%%%%%%

In the problem (\ref{ur1})--(\ref{ur3}) we replace the boundary
conditions (\ref{ur2}) with

\begin{equation}
\left\{
\begin{array}{rcl}
k(0,t)u_x(0,t)=\beta_1(t)u(0,t)-\mu_1(t),     \\
-k(l,t)u_x(l,t)=\beta_2(t)u(l,t)-\mu_2(t).%\quad 0\leq t\leq T.\\
\end{array}
\right. \label{ur11.1}
\end{equation}

In the rectangle $\bar Q_T$ we consider the Robin boundary value
problem (\ref{ur1}), (\ref{ur3}), (\ref{ur11.1}).

 {\bf Theorem 2.} If
$k(x,t)\in C^{1,0}(\bar Q_T)$, $q(x,t)$, $f(x,t)\in C(\bar Q_T)$,
$k(x,t)\geq c_1>0$, $q(x,t)\geq 0$ everywhere on $\bar Q_T$,
$\beta_i(t), \mu_i(t)\in C[0,T]$, $\beta_i(t)\geq \beta_0>0$, for
all $t\in [0,T]$, $i=1,2$, then the solution $u(x,t)$ of the problem
(\ref{ur1}), (\ref{ur3}), (\ref{ur11.1}) satisfies the a priori
estimate:

$$
\int\limits_{0}^{l}\Bbb{P}_{
(\omega)}^{(\theta-1)}\left(D_{0t}\right)u^2(x,t)dx+
\gamma_1\left(\int\limits_{0}^{t}\left(\|u_x(x,s)\|_0^2+
u^2(0,s)+u^2(l,s)\right)ds\right)\leq
$$
$$
\leq
\frac{\delta_1}{\gamma_1}\left(\int\limits_{0}^{t}\left(\|f(x,s)\|_0^2+\mu_1^2(s)+\mu_2^2(s)\right)ds\right)+
$$
\begin{equation}
+\int\limits_{0}^{l}u_0^2(x)dx\int\limits_{\alpha}^{\beta}d\gamma\sum\limits_{r=1}^{m}\frac{\omega_r(x,\gamma)t^{1-\theta_r(x,\gamma)}}{\Gamma(2-\theta_r(x,\gamma))},
 \label{ur11.2}
\end{equation}

where $\gamma_1=\min\{c_1,\beta_0\}$, $\delta_1=\max\{1+l,l^2\}$.

{\bf Proof.} Let us multiply the equation (\ref{ur1}) by $u(x,t)$
and integrate the resulting relation over $x$ from $0$ to $l$:
$$
\int\limits_{0}^{l}u(x,t)\Bbb{P}_{(\omega)}^{(\theta)}\left(\partial_{0t}\right)u(x,t)dx
-\int\limits_{0}^{l}u(x,t)(k(x,t)u_x(x,t))_xdx+
$$
\begin{equation}\label{ur11.3}
+\int\limits_{0}^{l}q(x,t)u^2(x,t)dx=\int\limits_{0}^{l}u(x,t)f(x,t)dx.
\end{equation}

Now we transform the terms of the identity (\ref{ur11.3}):
$$
\int\limits_{0}^{l}u(x,t)\Bbb{P}_{(\omega)}^{(\theta)}\left(\partial_{0t}\right)u(x,t)dx\geq
\frac{1}{2}\int\limits_{0}^{l}\Bbb{P}_{(\omega)}^{(\theta)}\left(\partial_{0t}\right)u^2(x,t)dx.
$$
$$
-\int\limits_0^lu(ku_x)_xdx= \beta_1(t)u^2(0,t)+\beta_2(t)u^2(l,t)-
\mu_1(t)u(0,t)-\mu_2(t)u(l,t)+\int\limits_0^lku_x^2dx,
$$
$$
\left|\int\limits_0^lufdx\right|\leq
\varepsilon\|u\|_0^2+\frac{1}{4\varepsilon}\|f\|_0^2,\quad
\varepsilon>0.
$$
From (\ref{ur11.3}), taking into account the transformations
performed, one arrives at the inequality

$$
\frac{1}{2}\int\limits_{0}^{l}\Bbb{P}_{(\omega)}^{(\theta)}\left(\partial_{0t}\right)u^2(x,t)dx+c_1\|u_x(x,t)\|_0^2
+\beta_0u^2(0,t)+\beta_0u^2(l,t)\leq
$$

\begin{equation}
\leq \mu_1(t)u(0,t)+\mu_2(t)u(l,t)+
\varepsilon\|u\|_0^2+\frac{1}{4\varepsilon}\|f\|_0^2.
 \label{ur11.4}
\end{equation}

Using the inequalities $\mu_1(t)u(0,t)\leq\varepsilon
u^2(0,t)+(4\varepsilon)^{-1}\mu_1^2(t)$,
$\mu_2(t)u(l,t)\leq\varepsilon
u^2(l,t)+(4\varepsilon)^{-1}\mu_2^2(t)$, $\varepsilon>0$;
$\|u(x,t)\|_0^2\leq l^2\|u_x(x,t)\|_0^2+l(u^2(0,t)+u^2(l,t))$ with
$\varepsilon={\gamma_1}/({2\delta_1})$, from (\ref{ur11.4}) one has
the following inequality
$$
\int\limits_{0}^{l}\Bbb{P}_{(\omega)}^{(\theta)}\left(\partial_{0t}\right)u^2(x,t)dx+\gamma_1\left(\|u_x(x,t)\|_0^2
+u^2(0,t)+u^2(l,t)\right)\leq
$$
\begin{equation}
\leq
\frac{\delta_1}{\gamma_1}\left(\|f(x,t)\|_0^2+\mu_1^2(t)+\mu_2^2(t)\right).
\label{ur11.5}
\end{equation}

Changing the variable $t$ by $s$ in inequality (\ref{ur11.5}) and
integrating it over $s$ from $0$ to $t$, we obtain the a priori
estimate (\ref{ur11.2}).

The uniqueness and the continuous dependence of the solution of
problem (\ref{ur1}), (\ref{ur3}), (\ref{ur11.1}) on the input data
follow from the a priori estimate (\ref{ur11.2}).

%%%%%%%%%%%%%%%%%%%%%%%%%%%%%%%%%%%%%%%%%%%%%%%%%%%%%%%%%%
\section{Boundary value problems in difference setting}
%%%%%%%%%%%%%%%%%%%%%%%%%%%%%%%%%%%%%%%%%%%%%%%%%%%%%%%%%

Suppose that a solution $u(x,t)\in C^{4,3}( Q_T)$ of the problem
(\ref{ur1})--(\ref{ur3}) exists, and the coefficients of the
equation (\ref{ur1}) and the functions $f(x,t)$, $u_0(x)$ satisfy
the smoothness conditions, required for the construction of
difference schemes with the order of approximation
$O(\tau^{2-\theta_{\max}}+h^2)$, where
$\theta_{\max}=\max\limits_{r,x,\gamma}\theta_r(x,\gamma)$.

In the rectangle $\bar Q_T$ we introduce the grid
$\bar\omega_{h\tau}=\bar\omega_{h}\times\bar\omega_{\tau}$, where
$\bar\omega_{h}=\{x_i=ih, i=0,1,\ldots,N, hN=l\}$,
$\bar\omega_{\tau}=\{t_j=j\tau, j=0,1,\ldots,j_0, \tau j_0=T\}$.

Before to turn to the approximation of the problem
(\ref{ur1})--(\ref{ur3}), let us find the discrete analog of the
$\Bbb{P}_{(\omega)}^{(\theta)}\left(\partial_{0t}\right)u(x,t)$. For
any fixed $\gamma\in[\alpha,\beta]$ and $r\in\{1,2,\ldots,m\}$, the
following equality takes place

$$
\partial_{0t_{j+1}}^{\theta_r(x_i,\gamma)}u(x_i,\eta)=
\frac{1}{\Gamma(1-\theta_r(x_i,\gamma))}\int\limits_{0}^{t_{j+1}}\frac{\frac{\partial}{\partial\eta}u(x_i,\eta)d\eta}{(t_{j+1}-\eta)^{\theta_r(x_i,\gamma)}}=
$$
$$
=\frac{1}{\Gamma(1-\theta_r(x_i,\gamma))}\sum\limits_{s=0}^{j}\int\limits_{t_s}^{t_{s+1}}\frac{\frac{\partial}{\partial\eta}u(x_i,\eta)d\eta}{(t_{j+1}-\eta)^{\theta_r(x_i,\gamma)}}=
$$
$$
=\frac{1}{\Gamma(1-\theta_r(x_i,\gamma))}\sum\limits_{s=0}^{j}\int\limits_{t_s}^{t_{s+1}}\frac{\frac{\partial}{\partial\eta}u(x_i,\eta)|_{\eta=t_{s+1/2}}}{(t_{j+1}-\eta)^{\theta_r(x_i,\gamma)}}d\eta+
$$
$$
=\frac{1}{\Gamma(1-\theta_r(x_i,\gamma))}\sum\limits_{s=0}^{j}\int\limits_{t_s}^{t_{s+1}}\frac{\frac{\partial^2}{\partial\eta^2}u(x_i,\eta)|_{\eta=t_{s+1/2}}\left(\eta-t_{s+1/2}\right)}{(t_{j+1}-\eta)^{\theta_r(x_i,\gamma)}}d\eta+O(\tau^2)=
$$
$$
=\frac{1}{\Gamma(2-\theta_r(x_i,\gamma))}\sum\limits_{s=0}^{j}\left(t_{j-s+1}^{1-\theta_r(x_i,\gamma)}-t_{j-s}^{1-\theta_r(x_i\gamma)}\right)\frac{u(x_i,t_{s+1})-u(x_i,t_s)}{\tau}+
$$
\begin{equation}
+\frac{1}{\Gamma(1-\theta_r(x_i,\gamma))}\sum\limits_{s=0}^{j}\frac{\partial^2}{\partial\eta^2}u(x_i,\eta)|_{\eta=t_{s+1/2}}\int\limits_{t_s}^{t_{s+1}}\frac{\left(\eta-t_{s+1/2}\right)d\eta}{(t_{j+1}-\eta)^{\theta_r(x_i,\gamma)}}+O(\tau^2).
\label{ur11.6}
\end{equation}

Since
$$
\left|\frac{1}{\Gamma(1-\theta_r(x_i,\gamma))}\sum\limits_{s=0}^{j}\frac{\partial^2}{\partial\eta^2}u(x_i,\eta)|_{\eta=t_{s+1/2}}\int\limits_{t_s}^{t_{s+1}}\frac{\left(\eta-t_{s+1/2}\right)d\eta}{(t_{j+1}-\eta)^{\theta_r(x_i,\gamma)}}\right|\leq
$$
$$
\leq\frac{M}{\Gamma(1-\theta_r(x_i,\gamma))}\sum\limits_{s=0}^{j}\left|\int\limits_{t_s}^{t_{s+1}}\frac{\left(\eta-t_{s+1/2}\right)d\eta}{(t_{j+1}-\eta)^{\theta_r(x_i,\gamma)}}\right|=
$$
$$
=\frac{M}{\Gamma(1-\theta_r(x_i,\gamma))}\sum\limits_{s=0}^{j}\left|\int\limits_{t_{s+1/2}}^{t_{s+1}}\frac{\left(\eta-t_{s+1/2}\right)d\eta}{(t_{j+1}-\eta)^{\theta_r(x_i,\gamma)}}
-\int\limits_{t_s}^{t_{s+1/2}}\frac{\left(t_{s+1/2}-\eta\right)d\eta}{(t_{j+1}-\eta)^{\theta_r(x_i,\gamma)}}\right|=
$$
$$
=\frac{2^{\theta_r(x_i,\gamma)}M\tau^{2-\theta_r(x_i,\gamma)}}{4\Gamma(1-\theta(x_i,\gamma))}\sum\limits_{s=0}^{j}\left(\int\limits_{0}^{1}\frac{zdz}{(2s+1-z)^{\theta_r(x_i,\gamma)}}
-\int\limits_{0}^{1}\frac{zdz}{(2s+1+z)^{\theta_r(x_i,\gamma)}}\right)=
$$
$$
=\frac{2^{\theta_r(x_i,\gamma)}M\tau^{2-\theta_r(x_i,\gamma)}}{4\Gamma(1-\theta_r(x_i,\gamma))}\int\limits_{0}^{1}z\sum\limits_{s=0}^{j}\left(\frac{1}{(2s+1-z)^{\theta_r(x_i,\gamma)}}-\frac{1}{(2s+1+z)^{\theta_r(x_i,\gamma)}}\right)dz=
$$
$$
=\frac{2^{\theta_r(x_i,\gamma)}M\tau^{2-\theta_r(x_i,\gamma)}}{4\Gamma(1-\theta_r(x_i,\gamma))}
\int\limits_{0}^{1}z\left(\frac{1}{(1-z)^{\theta_r(x_i,\gamma)}}-\frac{1}{(2j+1+z)^{\theta_r(x_i,\gamma)}}\right)dz-
$$
$$
-\frac{2^{\theta_r(x_i,\gamma)}M\tau^{2-\theta_r(x_i,\gamma)}}{4\Gamma(1-\theta_r(x_i,\gamma))}
\int\limits_{0}^{1}z\sum\limits_{s=1}^{j}\left(\frac{1}{(2s-1+z)^{\theta_r(x_i,\gamma)}}-\frac{1}{(2s+1-z)^{\theta_r(x_i,\gamma)}}\right)dz\leq
$$
$$
\leq\frac{2^{\theta_r(x_i,\gamma)}M\tau^{2-\theta_r(x_i,\gamma)}}{4\Gamma(1-\theta_r(x_i,\gamma))}
\int\limits_{0}^{1}\frac{zdz}{(1-z)^{\theta_r(x_i,\gamma)}}=
\frac{2^{\theta_r(x_i,\gamma)}M\tau^{2-\theta_r(x_i,\gamma)}}{4\Gamma(3-\theta_r(x_i,\gamma))}
\leq\frac{M\tau^{2-\theta_{\max}}}{2}
$$
with $M=\max\limits_{(x,t)\in\bar Q_T}|\frac{\partial^2}{\partial
t^2}u(x,t)|$, then multiplying (\ref{ur11.6}) by
$\omega_r(x_i,\gamma)$ and summing
 the result over $r$ from $1$ to $m$, then integrating over $\gamma$ from $\alpha$ to $\beta$ one finds
\begin{equation}
\Bbb{P}_{(\omega)}^{(\theta)}\left(\partial_{0t_{j+1}}\right)u(x_i,t)=
\Bbb{P}_{(\omega)}^{(\theta)}\left(\Delta_{0t_{j+1}}\right)u(x_i,t)+O(\tau^{2-\theta_{\max}}),
\label{ur11.65}
\end{equation}
where
$$\Bbb{P}_{
(\omega)}^{(\theta)}\left(\Delta_{0t_{j+1}}\right)u(x_i,t)=\int\limits_{\alpha}^{\beta}d\gamma\sum\limits_{r=1}^{m}\omega_r(x_i,\gamma)\Delta_{0t_{j+1}}^{\theta_r(x_i,\gamma)}u(x_i,t),
$$
$$
\Delta_{0t_{j+1}}^{\theta_r(x_i,\gamma)}u(x_i,t)=
\frac{1}{\Gamma(2-\theta_r(x_i,\gamma))}\sum\limits_{s=0}^{j}\left(t_{j-s+1}^{1-\theta_r(x_i,\gamma)}-t_{j-s}^{1-\theta_r(x_i\gamma)}\right)\frac{u(x_i,t_{s+1})-u(x_i,t_s)}{\tau}.
$$

{\bf Lemma 2.}  For any function $v(t)$ defined on the grid
$\bar\omega_\tau$ one has the inequalities
\begin{equation}
v^{j+1}\Bbb{P}_{
(\bar\omega)}^{(\bar\theta)}\left(\Delta_{0t_{j+1}}\right)v\geq
\frac{1}{2}\Bbb{P}_{
(\bar\omega)}^{(\bar\theta)}\left(\Delta_{0t_{j+1}}\right)v^2.
\label{ur11.7}
\end{equation}

{\bf Proof.}   For any fixed $\gamma\in[\alpha,\beta]$ and
$r\in\{1,2,\ldots,m\}$, relying on lemma 2~\cite{Alikh:12} one finds
the following inequality
\begin{equation}
v^{j+1}\Delta_{0t_{j+1}}^{\bar\theta_r(\gamma)}v\geq
\frac{1}{2}\Delta_{0t_{j+1}}^{\bar\theta_r(\gamma)}v^2+\frac{\tau^2\Gamma(2-\bar\theta_r(\gamma))}{2}\left(\Delta_{0t_{j+1}}^{\bar\theta_r(\gamma)}v\right)^2\geq
\frac{1}{2}\Delta_{0t_{j+1}}^{\bar\theta_r(\gamma)}v^2.
\label{ur11.8}
\end{equation}
 Multiplying (\ref{ur11.8}) by $\omega_1(\gamma)$ and summing
 the result over $r$ from $1$ to $m$, then integrating over $\gamma$ from $\alpha$ to $\beta$ one
arrives at (\ref{ur11.7}). The proof of the lemma 2 is complete.

%%%%%%%%%%%%%%%%%%%%%%%%%%%%%%%%%%%%%%%%%%%%%%%%%%%%%%%%%%%
\subsection{The Dirichlet boundary value problem}
%%%%%%%%%%%%%%%%%%%%%%%%%%%%%%%%%%%%%%%%%%%%%%%%%%%%%%%%%%

To problem (\ref{ur1})--(\ref{ur3}), we assign the difference
scheme:
\begin{equation}\label{ur12}
\Bbb{P}_{
(\omega)}^{(\theta)}\left(\Delta_{0t_{j+1}}\right)y_i=\Lambda
y^{j+1}_i+\varphi_i^{j+1},\quad i=1,2,\ldots,N-1, \quad
j=0,1,\ldots,j_0-1,
\end{equation}
\begin{equation}\label{ur13}
y(0,t)=0, \quad y(l,t)=0,\quad j=0,1,\ldots,j_0,
\end{equation}
\begin{equation}\label{ur14}
y(x,0)=u_0(x), \quad i=0,1,\ldots,N,
\end{equation}
where $\Lambda y=(ay_{\bar x})_x-dy$, $v_{\bar
x,i}=(v_{i}-v_{i-1})/h$, $v_{x,i}=(v_{i+1}-v_i)/h$,
$a_i^{j+1}=k(x_{i-1/2},t_{j+1})$, $d_i^{j+1}=q(x_i,t_{j+1})$,
$\varphi_i^{j+1}=f(x_i,t_{j+1})$,
$\Delta_{0t_j}^{\theta_r(x_i,\gamma)}y_i=
\sum\limits_{s=0}^{j}(t_{j-s+1}^{1-\theta_r(x_i,\gamma)}-t_{j-s}^{1-\theta_r(x_i,\gamma)})y_{t,i}^s/\Gamma(2-\theta_r(x_i,\gamma))$
is the difference analogue of the Caputo fractional derivative of
order $\theta_r(x_i,\gamma)$, $y_{t,i}^s=(y_i^{s+1}-y_i^s)/\tau$.

 According to \cite{Samar:77} and the formula (\ref{ur11.65}),
the order of the approximation of the difference scheme
(\ref{ur12})--(\ref{ur14}) is $O(\tau^{2-\theta_{\max}}+h^2)$.

{\bf Theorem 3.} The difference scheme (\ref{ur12})--(\ref{ur14}) is
absolutely stable and its solution satisfies the following a priori
estimate:
\begin{equation}\nonumber
\int\limits_{\alpha}^{\beta}d\gamma\sum\limits_{r=1}^{m}\left(\frac{\omega_r(x_i,\gamma)}{\Gamma(2-\theta_r(x_i,\gamma))},
\sum\limits_{s=0}^{j}(t_{j-s+1}^{1-\theta_r(x_i,\gamma)}-t_{j-s}^{1-\theta_r(x_i,\gamma)})(y_i^{s+1})^2\right)+
\end{equation}
\begin{equation}\label{ur15}
+c_1\sum\limits_{s=0}^{j}\|y_{\bar
x}^{j+1}]|_0^2\tau\leq\frac{l^2}{2c_1}\sum\limits_{s=0}^{j}\|\varphi^s\|_0^2\tau+\int\limits_{\alpha}^{\beta}d\gamma\sum\limits_{r=1}^{m}\left(
\frac{\omega_r(x_i,\gamma)t_{j+1}^{1-\theta_r(x_i,\gamma)}}{\Gamma(2-\theta_r(x_i,\gamma))},u_0^2(x_i)\right),
\end{equation}
where $(y,v)=\sum\limits_{i=1}^{N-1}y_iv_ih$,
$(y,v]=\sum\limits_{i=1}^{N}y_iv_ih$, $\|y\|_0^2=(y,y)$,
$\|y]|_0^2=(y,y]$.

{\bf Proof.} Let us multiply scalarly equation (\ref{ur12}) by
$y_i^{j+1}$:
\begin{equation}\label{ur16}
\left(y^{j+1},\Bbb{P}_{
(\omega)}^{(\theta)}\left(\Delta_{0t_{j+1}}\right)y_i\right)-(\Lambda
y^{j+1},y^{j+1})=(\varphi^{j},y^{j+1}).
\end{equation}

Let us transform the terms in identity (\ref{ur16}):
\begin{equation}\label{ur17}
-(\Lambda y^{j+1},y^{j+1})=(a,(y_{\bar
x}^{j+1})^2]+(d,(y^{j+1})^2\geq c_1\|y_{\bar x}^{j+1}]|_0^2,
\end{equation}
\begin{equation}\label{ur18}
|(\varphi^{j+1},y^{j+1})|\leq
\varepsilon\|y^{j+1}\|_0^2+\frac{1}{4\varepsilon}\|\varphi^{j+1}\|_0^2\leq
\frac{\varepsilon l^2}{2}\|y_{\bar
x}^{j+1}]|_0^2+\frac{1}{4\varepsilon}\|\varphi^{j+1}\|_0^2,\quad
\varepsilon>0.
\end{equation}
 Relying on lemma 2 one has
\begin{equation}\label{ur19}
\left(y^{j+1},\Bbb{P}_{
(\omega)}^{(\theta)}\left(\Delta_{0t_{j+1}}\right)y_i\right)\geq
\frac{1}{2}\left(1,\Bbb{P}_{
(\omega)}^{(\theta)}\left(\Delta_{0t_{j+1}}\right)y_i^2\right).
\end{equation}

 From (\ref{ur16}) with taking into account (\ref{ur17}),
(\ref{ur18}) and (\ref{ur19}), it follows that
\begin{equation}\label{ur23}
\frac{1}{2}\left(1,\Bbb{P}_{
(\omega)}^{(\theta)}\left(\Delta_{0t_{j+1}}\right)y_i^2\right)+c_1\|y_{\bar
x}^{j+1}]|_0^2\leq \frac{\varepsilon l^2}{2}\|y_{\bar
x}^{j+1}]|_0^2+\frac{1}{4\varepsilon}\|\varphi\|_0^2.
\end{equation}

Multiplying the inequality (\ref{ur23}) at $\varepsilon=c_1/l^2$, by
$\tau$ and summing over $s$ from $0$ to $j$, one obtains the a
priori estimate (\ref{ur15}).

The stability and convergence of the difference scheme
(\ref{ur12})--(\ref{ur14}) follow from the a priori estimate
(\ref{ur15}).

Here the results are obtained for the homogeneous boundary
conditions  $u(0,t)=0$, $u(l,t)=0$. In the case of inhomogeneous
boundary conditions $u(0,t)=\mu_1(t)$, $u(l,t)=\mu_2(t)$ the
boundary conditions for the difference problem will have the
following form:
\begin{equation}
y(0,t)=\mu_1(t),\quad y(l,t)=\mu_2(t). \label{ur2.2.16}
\end{equation}

Convergence of the difference scheme (\ref{ur12}), (\ref{ur14}),
(\ref{ur2.2.16}) follows from the results obtained above. Actually,
let us introduce the notation $y=z+u$. Then the error  $z=y-u$ is a
solution of the following problem:

 \begin{equation}\label{ur2.2.17}
\Bbb{P}_{
(\omega)}^{(\theta)}\left(\Delta_{0t_{j+1}}\right)z_i=\Lambda
z^{j+1} +\psi^{j+1},\,\, i=1,...,N-1,\,\, j=0,1,...,j_0-1,
\end{equation}

\begin{equation}
z(0,t)=0,\quad z(l,t)=0,\quad j=0,...,j_0, \label{ur2.2.18}
\end{equation}

\begin{equation}
z(x,0)=0,\quad i=0,...,N, \label{ur2.2.19}
\end{equation}
where $\psi\equiv\Lambda u^{j+1}-\Bbb{P}_{
(\omega)}^{(\theta)}\left(\Delta_{0t_{j+1}}\right)u_i+\varphi^{j+1}=O(\tau^{2-\theta_{\max}}+h^2)$.

The solution of the problem (\ref{ur2.2.17})--(\ref{ur2.2.19})
satisfies the estimation (\ref{ur15}) so that the solution of the
difference scheme (\ref{ur12}), (\ref{ur14}), (\ref{ur2.2.16})
converges to the solution of the corresponding differential problem
with order $O(\tau^{2-\theta_{\max}}+h^2)$.

%%%%%%%%%%%%%%%%%%%%%%%%%%%%%
\subsection{Numerical results}
%%%%%%%%%%%%%%%%%%%%%%%%%%%%%

 In this section, the following multi-term variable-distributed order time
fractional diffusion equation is considered:

\begin{equation}
\begin{cases}
\Bbb{P}_{
(\omega)}^{(\theta)}\left(\partial_{0t}\right)u(x,t)=\frac{\partial}{\partial
x}\left(k(x,t)\frac{\partial u}{\partial x}\right)-q(x,t)u+f(x,t),
\\
u(0,t)=\mu_1(t),\quad u(l,t)=\mu_2(t),\quad 0\leq t\leq 1,\\
u(x,0)=u_0(x),\quad 0\leq x\leq l,
\end{cases}
 \label{ur111}
\end{equation}
where
$$
\Bbb{P}_{ (\omega)}^{(\theta)}\left(\partial_{0t}\right)u(x,t)=
\int\limits_{0}^{1}d\gamma\sum\limits_{r=1}^{5}\omega_r(x,\gamma)\partial_{0t}^{\theta_r(x,\gamma)}u(x,t),
$$
$\theta_r(x,\gamma)=\frac{1+(rx+1)\gamma-\cos(rx\gamma)}{r+4}$,\quad
$\theta_{\max}\approx 0.856$,\quad
$\omega_r(x,\gamma)=(rx+1+rx\sin(rx\gamma))\frac{\Gamma(3-\theta_r(x,z))}{2r+8}$,
\quad $k(x,t)=\frac{8+\sin(t)}{3x^2+1}$,\quad
$q(x,t)=1-\sin{(xt)}$,\quad
$f(x,t)=\sum\limits_{r=1}^{5}\left(t^{2-\theta_r(x,0)}-t^{2-\theta_r(x,1)}\right)\frac{x^3+x+1}{\ln{t}}+(x^3+x+1)(t^2+1)(1-\sin(xt))$,\quad
$\mu_1(t)=t^2+1$,\quad $\mu_2(t)=3(t^2+1)$,\quad $u_0(x)=x^3+x+1$.

The exact solution is $u(x,t)=(x^3+x+1)(t^2+1)$.

A comparison of the numerical solution and exact solution is
provided in {\bf Table 1}.

{\bf Table 2} shows that when we take a fixed value $h=0.001$, then
as the number of time of our approximate scheme is decreased, a
redaction in the maximum error is observed, as expected and the
convergence order of time is $O(\tau^{2-\theta_{\max}})\approx
O(\tau^{1.144})$, where the convergence order is calculated by the
following formula: Convergence order $
=\log_{\frac{\tau_1}{\tau_2}}\frac{e_1}{e_2}$.

{\bf Table 3} shows that when we take
$h^2=\tau^{2-\theta_{\max}}\approx\tau^{1.144}$, as the number as
spatial subintervals/time steps is decreased, a reduction in the
maximum error is observed, es expected the convergence order of the
approximate scheme is $O(h^2)$, where the convergence order is
calculated by the following formula: Convergence order $
=\log_{\frac{h_1}{h_2}}\frac{e_1}{e_2}$.

\vspace{10mm}

\begin{tabular}{lc}
{\bf Table 1}\\
The error, numerical solution and exact solution, when $t=0.99$, $h=0.1$, $\tau=0.01$.\\
\hline
\hspace{1mm} Space($x_i$) \hspace{15mm}{Numerical solution} \hspace{15mm}{Exact solution} \hspace{15mm}{Error}\\
\hline
\hspace{2mm} 0.0000 \hspace{19mm} 1.9801000 \hspace{30mm}  1.9801000 \hspace{22mm}    0.0000000  \\
\hspace{2mm} 0.1000 \hspace{19mm} 2.1798574 \hspace{30mm}  2.1800901 \hspace{22mm}    0.0002337  \\
\hspace{2mm} 0.2000 \hspace{19mm} 2.3915106 \hspace{30mm}  2.3919608 \hspace{22mm}    0.0004502  \\
\hspace{2mm} 0.3000 \hspace{19mm} 2.6269549 \hspace{30mm}  2.6275927 \hspace{22mm}    0.0006378  \\
\hspace{2mm} 0.4000 \hspace{19mm} 2.8980853 \hspace{30mm}  2.8988664 \hspace{22mm}    0.0007811  \\
\hspace{2mm} 0.5000 \hspace{19mm} 3.2167962 \hspace{30mm}  3.2176625 \hspace{22mm}    0.0008663  \\
\hspace{2mm} 0.6000 \hspace{19mm} 3.5949810 \hspace{30mm}  3.5958616 \hspace{22mm}    0.0008806  \\
\hspace{2mm} 0.7000 \hspace{19mm} 4.0445325 \hspace{30mm}  4.0453443 \hspace{22mm}    0.0008118  \\
\hspace{2mm} 0.8000 \hspace{19mm} 4.5773422 \hspace{30mm}  4.5779912 \hspace{22mm}    0.0006490  \\
\hspace{2mm} 0.9000 \hspace{19mm} 5.2053012 \hspace{30mm}  5.2056829 \hspace{22mm}    0.0003817  \\
\hspace{2mm} 1.0000 \hspace{19mm} 5.9403000 \hspace{30mm}  5.9403000 \hspace{22mm}    0.0000000  \\
 \hline
\end{tabular}

\vspace{5mm}

\begin{tabular}{lc}
{\bf Table 2}\\
Maximum error behavior versus time grid size reduction at $t=0.99$ when $h=0.001$.\\
\hline
\hspace{2mm} $\tau$ \hspace{39mm}{Maximum error} \hspace{40mm}{Convergence order} \\
\hline
 0.99/10  \hspace{32mm} 0.0006796 \hspace{49mm}          \\
 0.99/20  \hspace{32mm} 0.0002474 \hspace{49mm}  1.458   \\
 0.99/40  \hspace{32mm} 0.0000907 \hspace{49mm}  1.448   \\
 \hline
\end{tabular}

\vspace{5mm}

\begin{tabular}{lc}
{\bf Table 3}\\
Maximum error behavior versus grid size reduction at $t=0.99$ when $h^2=\tau^{1.144}$.\\
\hline
\hspace{2mm} $h$ \hspace{39mm}{Maximum error} \hspace{40mm}{Convergence order} \\
\hline
\hspace{2mm} 1/10  \hspace{32mm} 0.0008536 \hspace{49mm}         \\
\hspace{2mm} 1/20  \hspace{32mm} 0.0002178 \hspace{49mm}  1.970  \\
\hspace{2mm} 1/40  \hspace{32mm} 0.0000550 \hspace{49mm}  1.985  \\
\hspace{2mm} 1/80  \hspace{32mm} 0.0000138 \hspace{49mm}  1.995  \\
 \hline
\end{tabular}

%%%%%%%%%%%%%%%%%%%%%%%%%%%%%%%%%%%%%%%%%%%%%%%%
\subsection{The Robin boundary value problem.}
%%%%%%%%%%%%%%%%%%%%%%%%%%%%%%%%%%%%%%%%%%%%%%%%

To the differential problem (\ref{ur1}), (\ref{ur3}), (\ref{ur11.1})
we assign the following difference scheme:
\begin{equation}\label{ur4.1}
\Bbb{P}_{
(\omega)}^{(\theta)}\left(\Delta_{0t_{j+1}}\right)y_i=\Lambda
y^{j+1}_i +\varphi^{j+1}, \quad i=0,\ldots,N, j=0,1,\ldots,j_0-1,
\end{equation}

\begin{equation}
y(x,0)=u_0(x),\quad i=0,\ldots,N, \label{ur4.2}
\end{equation}
where $\Lambda y=(a_1y_{x}-\tilde\beta_1 y)/(0.5h), i=0$, \quad
$\Lambda y=(ay_{\bar x})_x-dy, i=1,\ldots,N-1$, \quad $\Lambda
y=(-a_Ny_{\bar x}-\tilde\beta_2 y)/(0.5h), i=N$, \quad
$\varphi_0=(2\tilde\mu_1)/h$, \quad
$\varphi_N=(2\tilde\mu_2)/h$,\quad
$\tilde\beta_1=\beta_1+0.5hd_0$,\quad
$\tilde\beta_2=\beta_2+0.5hd_N$,\quad $\tilde\mu_1=\mu_1+0.5hf_0$,
\quad $\tilde\mu_2=\mu_2+0.5hf_N$.\\ The difference scheme
(\ref{ur4.1})--(\ref{ur4.2}) has the order of approximation
$O(\tau^{2-\theta_{\max}}+h^2)$.

{\bf Theorem 4.} The difference scheme (\ref{ur4.1})--(\ref{ur4.2})
  is absolutely stable and
its solution satisfies the following a priori estimate:

$$
\int\limits_{\alpha}^{\beta}d\gamma\sum\limits_{r=1}^{m}\left[\frac{\omega_r(x_i,\gamma)}{\Gamma(2-\theta_r(x_i,\gamma))},
\sum\limits_{s=0}^{j}(t_{j-s+1}^{1-\theta_r(x_i,\gamma)}-t_{j-s}^{1-\theta_r(x_i,\gamma)})(y_i^{s+1})^2\right]+
$$

$$
+ \gamma_1\sum\limits_{s=0}^{j}\left(\|(y_{\bar
x}^{s+1}]|_0^2+(y_0^{s+1})^2+(y_N^{s+1})^2\right)\tau\leq
$$

$$
\leq \frac{\delta_1}{\gamma_1}\sum\limits_{s=0}^{j}
\left((\tilde\mu_1^{s+1})^2+(\tilde\mu_2^{s+1})^2+\|\varphi^{s+1}\|_0^2\right)\tau+
$$

\begin{equation}\label{ur4.2.2}+
\int\limits_{\alpha}^{\beta}d\gamma\sum\limits_{r=1}^{m}\left[
\frac{\omega_r(x_i,\gamma)t_{j+1}^{1-\theta_r(x_i,\gamma)}}{\Gamma(2-\theta_r(x_i,\gamma))},u_0^2(x_i)\right],
\end{equation}

where $\gamma_1=\min\{c_1,\beta_0\}$, $\delta_1=\max\{1+l,l^2\}$,
$[y,v]=\sum_{i=1}^{N-1}y_iv_ih+0.5y_0v_0h+0.5y_Nv_Nh$,
$|[y]|_0^2=[y,y]$.

{\bf Proof.} Let us multiply scalarly equation (\ref{ur4.1}) by
$y^{j+1}$:
\begin{equation}\label{ur4.3}
\left[y^{j+1},\Bbb{P}_{
(\omega)}^{(\theta)}\left(\Delta_{0t_{j+1}}\right)y_i\right]-[\Lambda
y^{j+1},y^{j+1}] =[\varphi^{j+1},y^{j+1}],
\end{equation}
Let us transform the terms occurring in identity (\ref{ur4.3}) as
$$
\left[y^{j+1},\Bbb{P}_{
(\omega)}^{(\theta)}\left(\Delta_{0t_{j+1}}\right)y_i\right]
\geq\frac{1}{2}\left[1, \Bbb{P}_{
(\omega)}^{(\theta)}\left(\Delta_{0t_{j+1}}\right)y_i^2\right],
$$
$$
-[\Lambda y^{j+1},y^{j+1}]=\tilde\beta_1 (y_0^{j+1})^2+\tilde\beta_2
(y_N^{j+1})^2+(a,(y_{\bar x}^{j+1})^2]+[d,(y^{j+1})^2],
$$
$$
|[\varphi,y^{j+1}]|\leq
\varepsilon\|y^{j+1}\|_0^2+\tilde\mu_1y_0^{j+1}+\tilde\mu_2y_N^{j+1}+\frac{1}{4\varepsilon}\|\varphi\|_0^2,
\quad\varepsilon>0.
$$
Taking into account the above performed transformations, from
identity (\ref{ur4.3}) one arrives at the inequality
$$
\frac{1}{2}\left[1, \Bbb{P}_{
(\omega)}^{(\theta)}\left(\Delta_{0t_{j+1}}\right)y_i^2\right]+c_1\|y_{\bar
x}^{j+1}]|_0^2+\beta_0((y_0^{j+1})^2+(y_N^{j+1})^2)\leq
$$
\begin{equation}\label{ur4.4}
\leq
\varepsilon(\|y^{(\sigma)}\|_0^2+(y_0^{j+1})^2+(y_N^{j+1})^2)+\frac{1}{4\varepsilon}(\tilde\mu_1^2+\tilde\mu_2^2+\|\varphi\|_0^2).
\end{equation}

From (\ref{ur4.4}) at $\varepsilon=\gamma_1/(2\delta_1)$, using that
$\|y\|_0^2\leq l^2\|y_{\bar x}]|_0^2+l(y_0^2+y_N^2)$, one has the
following inequality:
$$
\left[1, \Bbb{P}_{
(\omega)}^{(\theta)}\left(\Delta_{0t_{j+1}}\right)y_i^2\right]+\gamma_1(\|y_{\bar
x}^{j+1}]|_0^2+(y_0^{j+1})^2+(y_N^{j+1})^2)\leq
$$
\begin{equation}\label{ur4.5}\leq
\frac{\delta_1}{\gamma_1}(\tilde\mu_1^2+\tilde\mu_2^2+\|\varphi\|_0^2).
\end{equation}
Multiplying inequality (\ref{ur4.5}) by $\tau$ and summing over $s$
from $0$ to $j$, one obtains a priori estimate (\ref{ur4.2.2}). The
stability and convergence of the difference scheme
(\ref{ur4.1})--(\ref{ur4.2}) follow from the a priori estimate
(\ref{ur4.2.2}).

%%%%%%%%%%%%%%%%%%%%%%%%%%%%%
\subsection{Numerical results}
%%%%%%%%%%%%%%%%%%%%%%%%%%%%%

In this section, the following multi-term variable-distributed order
diffusion equation is considered:

\begin{equation}
\begin{cases}
\Bbb{P}_{
(\omega)}^{(\theta)}\left(\partial_{0t}\right)u(x,t)=\frac{\partial}{\partial
x}\left(k(x,t)\frac{\partial u}{\partial x}\right)-q(x,t)u+f(x,t),
\\
k(0,t)u_x(0,t)=\beta_1(t)u(0,t)-\mu_1(t),     \\
-k(1,t)u_x(1,t)=\beta_2(t)u(1,t)-\mu_2(t), \quad 0\leq t\leq 1,\\
u(x,0)=u_0(x),\quad 0\leq x\leq 1,
\end{cases}
 \label{ur11111}
\end{equation}
where
$$
\Bbb{P}_{ (\omega)}^{(\theta)}\left(\partial_{0t}\right)u(x,t)=
\int\limits_{-2}^{3}d\gamma\sum\limits_{r=1}^{9}\omega_r(x,\gamma)\partial_{0t}^{\theta_r(x,\gamma)}u(x,t),
$$
$\theta_r(x,\gamma)=\frac{3+\gamma+e^{x(\gamma-3)}}{rx+14}$,\quad
$\theta_{\max}=0.5$,\quad
$\omega_r(x,\gamma)=\frac{(1+xe^{x(\gamma-3)})\Gamma(4-\theta_r(x,z))}{6rx+84}$,
\quad $k(x,t)=\frac{10+\cos(2t)}{5x^4+1}$,\quad
$q(x,t)=1-\cos{(2xt)}$,\quad
$f(x,t)=\sum\limits_{r=1}^{9}\left(t^{2-\theta_r(x,-2)}-t^{2-\theta_r(x,3)}\right)\frac{x^5+x+1}{\ln{t}}+(x^5+x+1)(t^3+1)(1-\cos(2xt))$,\quad
$\beta_1(t)=5+\cos(2t)$,\quad $\beta_2(t)=1-\cos(2t)/3$,\quad
$\mu_1(t)=-5(t^2+1)$,\quad $\mu_2(t)=13(t^3+1)$,\quad
$u_0(x)=x^5+x+1$.

The exact solution is $u(x,t)=(x^5+x+1)(t^3+1)$.

A comparison of the numerical solution and exact solution is
provided in {\bf Table 4}.

{\bf Table 5} shows that when we take a fixed value $h=0.01$, then
as the number of time of our approximate scheme is decreased, a
redaction in the maximum error is observed, as expected and the
convergence order of time is $O(\tau^{2-\theta_{\max}})=
O(\tau^{1.5})$, where the convergence order is calculated by the
following formula: Convergence order $
=\log_{\frac{\tau_1}{\tau_2}}\frac{e_1}{e_2}$.

{\bf Table 6} shows that when we take $h^2=\tau^{1.5}$, as the
number as spatial subintervals/time steps is decreased, a reduction
in the maximum error is observed, es expected the convergence order
of the approximate scheme is $O(h^2+\tau^{1.5})=O(h^2)$, where the
convergence order is calculated by the following formula:
Convergence order $ =\log_{\frac{h_1}{h_2}}\frac{e_1}{e_2}$.

\vspace{10mm}

\begin{tabular}{lc}
{\bf Table 4}\\
The error, numerical solution and exact solution, when $t=0.99$, $h=0.1$, $\tau=0.045$.\\
\hline
\hspace{1mm} Space($x_i$) \hspace{15mm}{Numerical solution} \hspace{15mm}{Exact solution} \hspace{15mm}{Error}\\
\hline
\hspace{2mm} 0.0000 \hspace{21mm} 1.9743424 \hspace{28mm}  1.9702990 \hspace{22mm}    0.0040434  \\
\hspace{2mm} 0.1000 \hspace{21mm} 2.1715724 \hspace{28mm}  2.1673486 \hspace{22mm}    0.0042238  \\
\hspace{2mm} 0.2000 \hspace{21mm} 2.3692967 \hspace{28mm}  2.3649893 \hspace{22mm}    0.0043074  \\
\hspace{2mm} 0.3000 \hspace{21mm} 2.5703757 \hspace{28mm}  2.5661765 \hspace{22mm}    0.0041992  \\
\hspace{2mm} 0.4000 \hspace{21mm} 2.7824031 \hspace{28mm}  2.7785945 \hspace{22mm}    0.0038086  \\
\hspace{2mm} 0.5000 \hspace{21mm} 3.0200717 \hspace{28mm}  3.0170203 \hspace{22mm}    0.0030514  \\
\hspace{2mm} 0.6000 \hspace{21mm} 3.3075397 \hspace{28mm}  3.3056889 \hspace{22mm}    0.0018508  \\
\hspace{2mm} 0.7000 \hspace{21mm} 3.6807935 \hspace{28mm}  3.6806565 \hspace{22mm}    0.0001370  \\
\hspace{2mm} 0.8000 \hspace{21mm} 4.1900054 \hspace{28mm}  4.1921658 \hspace{22mm}    0.0021604  \\
\hspace{2mm} 0.9000 \hspace{21mm} 4.9018802 \hspace{28mm}  4.9070100 \hspace{22mm}    0.0051298  \\
\hspace{2mm} 1.0000 \hspace{21mm} 5.9019826 \hspace{28mm}  5.9108970 \hspace{22mm}    0.0089144  \\
 \hline
\end{tabular}

\vspace{5mm}

\begin{tabular}{lc}
{\bf Table 5}\\
Maximum error behavior versus time grid size reduction at $t=0.99$ when $h=0.002$.\\
\hline
\hspace{2mm} $\tau$ \hspace{39mm}{Maximum error} \hspace{40mm}{Convergence order} \\
\hline
\hspace{2mm} 0.99/10 \hspace{32mm} 0.0051283 \hspace{49mm}         \\
\hspace{2mm} 0.99/20 \hspace{32mm} 0.0017203 \hspace{49mm}  1.576  \\
\hspace{2mm} 0.99/40 \hspace{32mm} 0.0005658 \hspace{49mm}  1.604  \\
\hspace{2mm} 0.99/60 \hspace{32mm} 0.0001834 \hspace{49mm}  1.625  \\
 \hline
\end{tabular}

\vspace{5mm}

\begin{tabular}{lc}
{\bf Table 6}\\
Maximum error behavior versus grid size reduction at $t=0.99$ when $h^2=\tau^{1.5}$.\\
\hline
\hspace{2mm} $h$ \hspace{39mm}{Maximum error} \hspace{40mm}{Convergence order} \\
\hline
\hspace{2mm} 1/10  \hspace{32mm} 0.0089152 \hspace{49mm}         \\
\hspace{2mm} 1/20  \hspace{32mm} 0.0022597 \hspace{49mm}  1.980  \\
\hspace{2mm} 1/40  \hspace{32mm} 0.0005727 \hspace{49mm}  1.980  \\
\hline
\end{tabular}

%%%%%%%%%%%%%%%%%%%%%%%%
\section{Conclusion}
%%%%%%%%%%%%%%%%%%%%%%%%%

%The results obtained in the present paper allow to apply the method
%of energy inequalities to finding a priori estimates for boundary
%value problems for the fractional diffusion equation in differential
%and difference settings exactly as in the classical case
%($\alpha(x)=1$). It is interesting to note that the condition
%$\sigma\geq 1/(3-2^{1-\alpha_i})$ at $\alpha(x)=1$ turns into the
%well known condition $\sigma\geq 1/2$ of the absolute stability of
%the difference schemes for the classical diffusion equation.
Solutions of the Dirichlet and Robin boundary value problems for the
multi-term variable-distributed order diffusion equation are
studied. A priori estimates for the corresponding differential and
difference problems are obtained by using the method of the energy
inequalities. The stability and convergence of the difference
schemes follow from these a priory estimates. The credibility of the
obtained results is verified by performing numerical calculations
for test problems.

The method of the energy inequalities proposed in the present paper
can be used to find a priori estimates for solutions of a wide class
of boundary value problems for the multi-term variable-distributed
order diffusion equation to which the maximum principle is not
applicable (for example, problems considered in
\cite{Alikh:08,Alikh:10_2,Alikh:13,Alikh:13_2}).

It should be emphasized that  1) from the considered equation at
$m=1$, $\theta_1(x,\gamma)=\gamma$,
$\omega_1(x,\gamma)=\omega(\gamma)$ one obtains the distributed
order time-fractional diffusion equation, 2) setting
$\theta_r(x,\gamma)=\theta_r=const$,
$\int_{\alpha}^{\beta}\omega_r(x,\gamma)d\gamma=\lambda_r(x)$ yields
the multi-term time-fractional diffusion equation.

%The present paper studies the equation with the variable-distributed
%order derivative
%generalizing the distributed order derivative \\
%$\int_{\alpha}^{\beta}\omega(\gamma)\partial_{0t}^\gamma ud\gamma$.
%In the presented case the order of the derivative depends on the
%function $\theta(x,\gamma)$, while $\alpha$ and $\beta$
%($\alpha<\beta$) may take arbitrary values including infinity. As an
%example, a test problem with the Robin boundary condition is solved
%for $\alpha=-2$, $\beta=3$. It should be emphasized that  1) at
%$\theta(x,\gamma)=\gamma$, $\omega(x,\gamma)=\omega(\gamma)$ one
%obtains the distributed order derivative, 2) setting
%$\theta(x,\gamma)=\theta \,\,\,(\theta=const)$,
%$\int_{\alpha}^{\beta}\omega(x,\gamma)d\gamma=1$ yields the
%fractional derivative of the order $\theta$.

\section{Acknowledgements}
This work was supported by the Russian Foundation for Basic Research
(project 14-01-31246).



%% References



%%%%%%%%%%%%%%%%%%%%%%%%%%%%%%%%%%%

\end{document}